\documentclass{siamart190516}

\usepackage{lipsum}
\usepackage{amsfonts}
\usepackage{graphicx}
\usepackage{epstopdf}
\usepackage{algorithmic}
\ifpdf
  \DeclareGraphicsExtensions{.eps,.pdf,.png,.jpg}
\else
  \DeclareGraphicsExtensions{.eps}
\fi

\newsiamremark{hypothesis}{Hypothesis}
\crefname{hypothesis}{Hypothesis}{Hypotheses}
\newsiamthm{claim}{Claim}

\headers{Least-squares neural network method}{Z. Cai, J. Chen AND M. Liu}

\title{Least-Squares ReLU Neural Network (LSNN) Method \\[1mm] for
Scalar Nonlinear Hyperbolic Conservation Law\thanks{This work was supported in part by the National Science Foundation under grant DMS-2110571.}
}

\author{Zhiqiang Cai\thanks{Department of Mathematics, Purdue University, 150 N. University Street, West Lafayette, IN 47907-2067 
  (\email{caiz@purdue.edu}, \email{chen2042@purdue.edu}).}
\and Jingshuang Chen\footnotemark[2]
\and Min Liu\thanks{School of Mechanical Engineering, Purdue University, 585 Purdue Mall,
West Lafayette, IN 47907-2088(\email{liu66@purdue.edu}). }}

\usepackage{amsopn}

\ifpdf
\hypersetup{
  pdftitle={FOSLS},
  pdfauthor={Z.Cai, J.Chen, M.Liu}
}
\fi

\usepackage[utf8]{inputenc}
\usepackage{bm}
\usepackage{mathtools}
\usepackage{indentfirst}
\usepackage{subfigure}
\usepackage{tcolorbox}
\usepackage{color}
\usepackage[export]{adjustbox}
\usepackage{scalerel}
\usepackage{multirow}
\usepackage{algorithmic}
\Crefname{ALC@unique}{Line}{Lines}

\newcommand{\R}{\mathbb{R}}

\usepackage{graphicx}
\usepackage{diagbox}
\usepackage{pifont}
\newcommand{\vertiii}[1]{{\left\vert\kern-0.25ex\left\vert\kern-0.25ex\left\vert #1 
    \right\vert\kern-0.25ex\right\vert\kern-0.25ex\right\vert}}

\newcommand{\btheta}{\mbox{\boldmath${\theta}$}}

\newcommand{\bomega}{\mbox{\boldmath$\omega$}}

\usepackage{enumitem}
\setlist[itemize]{left=16pt} 

\def\bb{{\bf b}}

\def\be{{\bf e}}
\def\bff{{\bf f}}

\def\bh{\small{\bf h}}
\def\bx{{\bf x}}

\def\bzz{{\bf z}}

\def\cM{{\cal M}}
\def\cN{{\cal N}}

\def\cS{{\cal S}}
\def\cT{{\cal T}}

\begin{document}

\maketitle

\begin{abstract}
In \cite{Cai2021linear}, we introduced the least-squares ReLU neural network (LSNN) method for solving the linear advection-reaction problem with discontinuous solution and showed that the method outperforms mesh-based numerical methods in terms of the number of degrees of freedom. This paper studies the LSNN method for scalar nonlinear hyperbolic conservation law. 
The method is a discretization of an equivalent least-squares (LS) formulation in the set of neural network functions with the ReLU activation function. Evaluation of the LS functional is done by using numerical integration and conservative finite volume scheme.
Numerical results of some test problems show that the method is capable of approximating the discontinuous interface of the underlying problem automatically through the {\it free} breaking lines of the ReLU neural network. Moreover, the method does not exhibit the common Gibbs phenomena along the discontinuous interface. 
\end{abstract}



\section{Introduction}

Let $\Omega$ be a bounded domain in ${\R}^d$ ($d=1, \,2$, or $3$) with Lipschitz boundary, 
consider the scalar nonlinear hyperbolic conservation law 
\begin{equation} \label{pde}
    \left\{\begin{array}{rcll}
    u_t(\bx,t) + \nabla_\bx \cdot \bff (u) &= & 0, &\text{ in }\,\, {\Omega}  \times I, \\[2mm]
    u&=&
   {g}, &\text{ on }\,\, {\Gamma}_{-} ,\\[2mm]
    u(\bx,0) &=& u_0(\bx), &\text{ in }\,\, {\Omega},
    \end{array}\right.
\end{equation}
where $u_t$ is the partial derivative of $u$ with respect to temporal variable $t$; $\nabla_\bx\cdot$ is a divergence operator with respect to spatial variable $\bx$; $\bff (u) = (f_1(u),...,f_d(u))$ is the spatial flux vector field; $I=(0, T)$ is temporal interval; ${\Gamma}_-$ is the part of the boundary $\partial {\Omega} \times I$ where the characteristic curves enter the domain ${\Omega}  \times I$; and the boundary data ${g}$ and the initial data $u_0$ are given scalar-valued functions. 

Numerical methods for (\ref{pde}) have been intensively studied during the past several decades by many researchers and many numerical schemes have been developed. A major difficulty in simulation is that the solution of a scalar hyperbolic conservation law is often discontinuous due to the discontinuous initial/boundary condition or shock formation; moreover, there is no {\it a priori} knowledge of the location of the discontinuities. It is well-known that 
traditional mesh-based numerical methods often exhibit oscillations near a discontinuity (called the Gibbs phenomena).
Such spurious oscillations are unacceptable for many applications (see, e.g, \cite{hesthaven2007nodal}). 
To eliminate or reduce the Gibbs phenomena, finite difference and finite volume methods often use numerical techniques such as limiters and filters and conservative schemes such as Roe, ENO/WENO, etc. have been developed \cite{gottlieb1997gibbs,hesthaven2017numerical,hesthaven2007nodal, leveque1992numerical}; and finite element methods 
usually employ discontinuous finite elements \cite{brezzi2004discontinuous,dahmen2012adaptive, demkowicz2010class} and/or adaptive mesh refinement (AMR) to generate locally refined elements along discontinuous interfaces (see, e.g., \cite{burman2009posteriori, houston1999posteriori, houston2000posteriori}). 

Recently, there has been increasing interests in using neural networks (NNs) to solve partial differential equations (see, e.g., \cite{PNAS2019, cai2020deep, Weinan18, raissi2019physics,Sirignano18}). NNs produce a large class of functions through compositions of linear transformations and activation functions. One of the striking features of NNs is that this class of functions is not subject to a hand-crafted geometric mesh or point cloud as are the traditional, well-studied finite difference, finite volume, and finite element methods. The physical partition of the domain $\Omega$, formed by free hyper-planes, can automatically adapt to the target function. 
To make use this powerful approximation property of NNs, in \cite{Cai2021linear}, we studied least-squares neural network (LSNN) method for solving linear advection-reaction problem with discontinuous solution. The LSNN method is based on a direct application of the lease-squares principle to the underlying problem
studied in (\cite{bochev2001improved, de2004least}) and on the ReLU neural network as the approximation class of functions. Compared to various AMR methods that locate the discontinuous interface through local mesh refinement, the LSNN method is much more effective in terms of the number of the degrees of freedom. 

The purpose of this paper is to study the space-time LSNN method for solving the scalar hyperbolic conservation law. For the nonlinear hyperbolic conservation law, the differential equation is not generally sufficient to determine the solution. An additional constraint the so-called Rankine-Hugoniot (RH) jump condition \cite{godlewski2013numerical,leveque1992numerical, thomas2013numerical}, is needed at where the solution is not continuous. To enforce this condition weakly, \cite{de2005numerical} introduced an independent variable, the spatial-temporal flux, for the inviscid Burgers equation and applied the least-squares principle to the resulting equivalent system. A variant of this method was also studied in \cite{ de2005numerical,manteuffel2020least} by using the Helmholtz decomposition of the flux. 

Due to the training difficulty of the least-squares method of \cite{de2005numerical}, in this paper we employ the naive least-squares method used for the linear advection-reaction problem \cite{Cai2021linear}, i.e., a direct application of least-squares principle to the PDE, initial and inflow boundary conditions. To ensure that the numerical solution enforces the RH jump condition, we introduce implicit discrete finite difference operators in section 3 by following ideas of the explicit conservative schemes such as Roe, ENO, etc.

Based on our numerical experience, it is difficult to train the network for problems with sharp changes when the LSNN method is applied to the entire computational domain $\Omega$.
This is due to the nature of nonlinear hyperbolic conservation laws since information is transported from initial and inflow boundary to the rest of the domain along the flow direction. To overcome this training difficulty, we propose the block space-time LSNN method in section 4. Basically, we partition the computational domain into a number of blocks based on the ``inflow'' boundary and initial conditions, then the method solves problems on these blocks sequentially. The trained parameters of the NN model for the previous block is used as an initial for the current block.

NN-based numerical methods for solving scalar nonlinear hyperbolic conservation laws have been studied recently by many researchers (see, e.g., \cite{PNAS2019, raissi2017physics}). The most popular one is the so-called the physics informed neural network (PINN) method proposed in \cite{raissi2017physics}. Basically, the PINN is a archaic version of the NN method of LS type; it employs a primitive form of least-squares formulation and uses the automatic differentiation method to differentiate the neural network. This is why the PINN and other NN-based methods are only applicable to problems with smooth solution such as the viscous Burgers equation but not the inviscid one studied in this paper.


The paper is organized as follows. The ReLU neural network and the space-time LSNN method are introduced in section 2. Conservative finite difference operators are discussed in section 3. The block LSNN method
is explained in section 4. Finally, implementation and numerical results for various one dimensional benchmark test problems are presented in section 5.

\section{Space-Time Least-Squares Neural Network Method}

In this section, we describe least-squares neural network method for the scalar hyperbolic conservation law. 

A deep neural network (DNN) defines a scalar-valued function
\[
\cN:\, \bzz=(\bx,t)\in\R^{d+1}
\longrightarrow \cN(\bzz)\in\R.
\]
A DNN function $\cN(\bzz)$ is typically represented as compositions of many layers of functions:  
\begin{equation}\label{DNN}
 \cN(\bzz)=\bomega^{(L)} \left(\mathcal{N}^{(L-1)} \circ \cdots \mathcal{N}^{(2)}\circ \mathcal{N}^{(1)}(\bzz)\right)- b^{(L)},
\end{equation}
where  $\bomega^{(L)}\in \R^{n_{L-1}}$, $b^{(L)}\in \R$, the symbol $\circ$ denotes the composition of functions, and $L$ is the depth of the network. For $l=1,\,\cdots,\,L-1$, the $\mathcal{N}^{(l)}: \R^{n_{l-1}} \rightarrow \R^{n_{l}}$ is called
the $l^{th}$ hidden layer of the network defined as follows:
\begin{equation}\label{layerdef}
  \mathcal{N}^{(l)}(\bzz^{(l-1)})
  = \sigma (\bomega^{(l)}\bzz^{(l-1)}-\bb^{(l)})
  \quad\mbox{for } \bzz^{(l-1)}\in \R^{n_{l-1}},
\end{equation}
where $\bomega^{(l)} 
\in \R^{n_{l}\times n_{l-1}}$, $\bb^{(l)}\in \R^{n_{l}}$, $\bzz^{(0)}=\bzz$, and $\sigma$ is the activation function and its application 
to a vector is defined component-wisely. This paper will use the popular rectified linear unit (ReLU) activation function defined by
\begin{equation}\label{tau-k}
 \sigma(s) = \max\{0,\,s\}
 =\left\{\begin{array}{rclll}
 0, & \mbox{if }  s\leq 0,\\[2mm]
 s, & \mbox{if } s >0.
 \end{array}\right.
 \end{equation}
Denote the set of NN functions by
\[
\cM({\small\btheta},L)=\big\{\cN(\bzz)=\bomega^{(L)} \left(\mathcal{N}^{(L-1)} \circ \cdots \mathcal{N}^{(2)}\circ \mathcal{N}^{(1)}(\bzz)\right)- b^{(L)} :\,  \bomega^{(l)} 
\in \R^{n_{l}\times n_{l-1}},\,\, \bb^{(l)}\in \R^{n_{l}} 
\big\},
\]
where $\mathcal{N}^{(l)}(\bzz^{(l-1)})$ is defined in (\ref{layerdef}) and ${\small\btheta}$ denotes all parameters: $\bomega^{(l)}$ and $\bb^{(l)}$ for $l=1,...,L$.
It is easy to see that $\cM({\small\btheta},L)$ is a set, but not a linear space. 

Applying the least-squares principle directly to the problem in (\ref{pde}), we have the following 
least-squares (LS) functional
 \begin{equation}\label{ls}
    \mathcal{L}(v;{ g}) = \| v_t + \nabla_\bx \cdot \bff (v)\|_{0,\Omega\times I}^2 +  \|v-g\|_{0, \Gamma_-}^2 + \|v(\bx,0)-u_0(\bx)\|_{0, \Omega}^2 .
\end{equation}
Then the least-squares approximation is to find $u_{_N}(\bx,t;{\small\btheta}^*) \in \cM({\small\btheta},L)$ such that 
\begin{equation}\label{L-NN}
     \mathcal{L}\big(u_{_N}(\cdot;{\small\btheta}^*);\,{\bf f}\big)
     = \min\limits_{v\in \cM({\scriptsize\btheta},L)} \mathcal{L}\big(v(\cdot;{\small\btheta});\,{g}\big)
     = \min_{{\scriptsize \btheta}\in\R^{N}}\mathcal{L}\big(v(\cdot;{\small\btheta});\,{g}\big),
\end{equation}
where $N$ is the total number of parameters in $\cM({\small\btheta},L)$ given by
\[
N=M_d(L) =\sum^L_{l=1} n_{l}\times (n_{l-1}+1).
\]

Similar to \cite{Cai2021linear,cai2020deep}, the integral in the LS functional is evaluated by numerical integration. To do so, let 
\[
{\cal T}=\{K :\, K\mbox{ is an open subdomain of } \Omega\times I \}
\] 
be a partition of the domain $\Omega$. Then
\[
{\cal E}_{-}=\{E
=\partial K \cap \Gamma_-:\,\, K\in\mathcal{T}\}
\quad\mbox{and}\quad
{\cal E}_0=\{E
=\partial K \cap \left(\Omega\times \{0\}\right):\,\, K\in\mathcal{T}\}
\]
are partitions of the boundary $\Gamma_-$ and $\Omega\times \{0\}$, respectively.
Let $\bzz_{_K}=(\bx_{_K}, t_{_K})$ and $\bzz_{_E}=(\bx_{_E}, t_{_E})$ be the centroids of $K\in {\cal T}$ and $E$ in ${\cal E}_-$ or ${\cal E}_0$, respectively. Define the discrete LS functional as follows:
\begin{equation}\label{L-NN-d}
\begin{split}
     \mathcal{L}_{_{\small {\cal T}}}\big(v(\cdot; {\small\btheta});{g}\big) 
     = \sum_{K \in {\cal T}}  \big( \delta_\tau v + \nabla_{\bh} \!\cdot \bff (v)\big)^2(\bzz_{_K}; {\small\btheta})\,|K|+  \sum_{E\in {\cal E}_-} \big(v-g\big)^2(\bzz_{_E}; {\small\btheta})|E| +  \sum_{E\in {\cal E}_0} \big(v-u_0\big)^2(\bzz_{_E}; {\small\btheta})|E|,
\end{split}
\end{equation}
where $|K|$ and $|E|$ are the $d$ and $d-1$ dimensional measures of $K$ and $E$, respectively; $\delta_\tau$ and $\nabla_{\bh} \!\cdot$ are finite difference operators to be defined in the subsequent section.
Then the discrete least-squares approximation is to find ${u}_{_{\small {\cal T}}}(\bzz,{\small\btheta}^*)\in \cM({\small\btheta},L)$ such that
 \begin{equation}\label{discrete_minimization_functional}
  \mathcal{L}_{_{\small {\cal T}}} \big({u}_{_{\small {\cal T}}}(\cdot,{\small\btheta}^*);{g}\big) 
  = \min\limits_{v\in \cM({\scriptsize\btheta},L)} \mathcal{L}_{_{\small {\cal T}}}\big(v(\cdot;{\small\btheta});\,{g}\big)
 = \min_{{\scriptsize \btheta}\in\R^{N}}\mathcal{L}_{_{\small {\cal T}}} \big(v(\cdot; {\small\btheta});{g}\big).
\end{equation}

\section{Conservative Finite Volume Operator}

How to discretize the differential operator is critical for the success of the LSNN method. Using finite difference quotient along coordinate directions to approximate the differential operator usually results in very poor numerical approximation. To overcome this difficulty, we employ  conservative finite volume schemes to evaluate the derivatives in the least-squares functional in (\ref{L-NN-d}). There are many conservative schemes
such as Roe's scheme, ENO, and WENO, etc. (see, e.g., \cite{roe1981approximate, shu1998essentially, shu1988efficient}). 
For simplicity, we briefly describe the finite volume operator using the idea of either Roe's scheme or second-order accurate ENO in this section. Note that the finite volume operators described in this section are implicit in time. 

For any $K\in\cT$, let $(\bx_{_K},t_{_K})$ be the centroid of $K$. Let $(\bh,\tau)=(h_1, ...,h_d,\tau)$ be a step size such that $(\bx_{_K}\pm\bh/2,t_{_K}\pm \tau/2)\in K$. For $i=1,..., d$, let $\bh_i = h_i \be_i$, where $\be_i$ is the unit vector in the $x_i$-coordinate direction. Then the
finite volume operator of the least-squares functional in (\ref{L-NN-d}) at the point $(\bx_{_K},t_{_K})$
is given by
\begin{eqnarray}\nonumber   
    && \quad  \big(\delta_\tau v+ \nabla_{\bh} \cdot \bff (v)\big)(\bx_{_K},t_{_K})\\[2mm] \label{Roe}
     &=&\dfrac{v\big(\bx_{_K},t_{_K}\big)-v\big(\bx_{_K},t_{_K}-\frac{1}{2}\tau\big)}{\tau} + \sum_{i=1}^{d}\dfrac{\hat{f}_i\big(v(\bx_{_K}+\frac{1}{2}\bh_i,t_{_K})\big)-\hat{f}_i\big(v(\bx_{_K}-\frac{1}{2}\bh_i,t_{_K})\big)}{h_i},
\end{eqnarray}
where $\hat{f}_i\big(v(\bx_{_K}\pm\frac{1}{2}\bh_i,t_{_K})\big)$ are the $i^{th}$ component of the numerical flux at $(\bx_{_K}\pm\frac{1}{2}\bh_i,t_{_K})$. Various conservative schemes are more or less on how to reconstruct proper numerical flux. 

Below, we describe the numerical fluxes by Roe and ENO. To this end, we introduce the Roe speed at point $(\bx_{_K}\pm\frac{1}{2}\bh_i,t_{_K})$ in the $\be_i$ direction 
\begin{equation}\label{roe_a}
 {a}_i\left(\bx_{_K}\!\pm\frac{1}{2}\bh_i,t_{_K}\right)=\left\{\begin{array}{ll}
 \dfrac{f_i\big(v(\bx_{_K}\pm\bh_i,t_{_K})\big)-f_i\big(v(\bx_{_K},t_{_K})\big) }{v\big(\bx_{_K}\pm\bh_i,t_{_K}\big) -v(\bx_{_K},t_{_K})}, & \text{if } v\big(\bx_{_K}\pm\bh_i,t_{_K}\big) \neq v(\bx_{_K},t_{_K}),\\[2mm]
 \noalign{\vskip9pt}
 f_i^\prime\big(v(\bx_{_K},t_{_K})\big), & \text{if }  v\big(\bx_{_K}\pm\bh_i,t_{_K}\big) = v(\bx_{_K},t_{_K}).
 \end{array}
 \right.
\end{equation}
Then the $i^{th}$ components of the Roe numerical flux at $(\bx_{_K}\pm\frac{1}{2}\bh_i,t_{_K})$ are given by
\begin{eqnarray}\nonumber
    && \quad \hat{f}_i\left(v\big(\bx_{_K}\pm\frac{1}{2}\bh_i,t_{_K}\big)\right)\\[2mm] \label{roe1}
    \qquad &=& \dfrac{{f}_i\big(v(\bx_{_K},t_{_K})\big)+{f}_i\big(v(\bx_{_K}\pm\bh_i,t_{_K})\big)}{2} \mp  \left|a_i\left(\bx_{_K}\pm\frac{1}{2}\bh_i,t_{_K}\right)\right|\dfrac{v\big(\bx_{_K}\pm\bh_i,t_{_K}\big) -v(\bx_{_K}\!,t_{_K})}{2}.
\end{eqnarray}

The key idea of the Roe's scheme is to use only grid points, if possible, on one side of the interface for constructing a finite volume scheme so that the RH condition is not violate. This is done through the signs of the Roe speed $a_i$ at midpoints. This key idea was further explored for developing higher order schemes, e.g., the ENO schemes introduced in \cite{harten1987uniformly} (see also \cite{shu1998essentially,shu1988efficient}), by employing extra grid points. To make sure all used grid points locate on one side of the interface, it requires additional decisions and, hence, the ENO schemes are generally sophisticated. 

For simplicity, we describe the second order ENO numerical flux here. The ENO uses the sign of the Roe speed to build up upwind scheme. Specifically, 
 \begin{equation}\label{upwind}
    \hat{f}_i\left(v(\bx_{_K}+\frac{1}{2}\bh_i,t_{_K})\right)=\left\{\begin{array}{ll}
 \hat{f}_i^{-}\left(v(\bx_{_K}+\frac{1}{2}\bh_i,t_{_K})\right), & \text{if } {a}_i\left(\bx_{_K}+\frac{1}{2}\bh_i,t_{_K}\right)\ge 0,\\[4mm]
 \hat{f}_i^{+}\left(v(\bx_{_K}+\frac{1}{2}\bh_i,t_{_K})\right), & \text{if } {a}_i\left(\bx_{_K}+\frac{1}{2}\bh_i,t_{_K}\right)<0.
 \end{array}\right. 
 \end{equation}

Additionally, the ENO uses the magnitudes of the finite difference quotient of the $i^{th}$ component of the flux with respect to $x_i$ over the neighboring intervals to determine which side of grid points are used. In this way, the ENO again tries to use grid points on one side of the discontinuity if possible. More precisely, let 
\begin{eqnarray*}
f_i(\bx_{_K},t_{_K};\bh_i)
&:=& \dfrac{f_i\big(v(\bx_{_K}+\bh_i,t_{_K})\big)-f_i\big(v(\bx_{_K},t_{_K})\big)}{h_i}\\[2mm]
\mbox{and }\,\, f_i(\bx_{_K},t_{_K};-\bh_i)
&:=& \dfrac{f_i\big(v(\bx_{_K},t_{_K})\big)-f_i\big(v(\bx_{_K}-\bh_i,t_{_K})\big)}{h_i}.
\end{eqnarray*}
In the case that ${a}_i\left(\bx_{_K}+\frac{1}{2}\bh_i,t_{_K}\right)\ge 0$, combining with (\ref{upwind}), the ENO numerical flux is given by
 \begin{eqnarray}\nonumber    
   \quad &&  \hat{f}_i^-\left(v(\bx_{_K}+\frac{1}{2}\bh_i,t_{_K})\right) \\[2mm] \label{n-flux}
   &=& \left\{\begin{array}{ll}
 -\frac{1}{2}f_i\big(v(\bx_{_K}-\bh_i,t_{_K})\big)   +\frac{3}{2} f_i\big(v(\bx_{_K},t_{_K})\big), & \text{if } \big|f_i(\bx_{_K},t_{_K};-\bh_i)\big| < \big|f_i(\bx_{_K},t_{_K};\bh_i)\big|,\\[4mm]
 \frac{1}{2}f_i\big(v(\bx_{_K},t_{_K})\big) + \frac{1}{2} f_i\big(v(\bx_{_K}+\bh_i,t_{_K})\big), & \text{if } \big|f_i(\bx_{_K},t_{_K};-\bh_i)\big| \ge \big|f_i(\bx_{_K},t_{_K};\bh_i)\big|.
 \end{array}\right. 
 \end{eqnarray}
If ${a}_i\left(\bx_{_K}+\frac{1}{2}\bh_i,t_{_K}\right)< 0$, then the numerical flux is reconstructed by
   \begin{eqnarray} \nonumber
    && \qquad \quad \hat{f}_i^+\left(v(\bx_{_K}+\frac{1}{2}\bh_i,t_{_K})\right) \\[2mm] \label{n-flux-1}
\!\!\! &=& \left\{\!\!\begin{array}{ll}\nonumber 
 \frac{1}{2}f_i\big(v(\bx_{_K},t_{_K})\big) + \frac{1}{2} f_i\big(v(\bx_{_K}\!\!+\bh_i,t_{_K})\big), & \! \text{if } \big|f_i(\bx_{_K}\!\!+\bh_i,t_{_K};-\bh_i)\big| < \big|f_i(\bx_{_K}\!\!+\bh_i,t_{_K};\bh_i)\big|,\\[4mm]
 \frac{3}{2}f_i\big(v(\bx_{_K}\!\!+\bh_i,t_{_K})\big) -\frac{1}{2} f_i\big(v(\bx_{_K}\!\!+2\bh_i,t_{_K})\big), &\! \text{if } \big|f_i(\bx_{_K}\!\!+\bh_i,t_{_K};-\bh_i)\big| \geq \big|f_i(\bx_{_K}\!\!+\bh_i,t_{_K};\bh_i)\big|.
 \end{array}\right. 
 \end{eqnarray}
In a similar fashion,  $\hat{f}_i\left(v(\bx_{_K}-\frac{1}{2}\bh_i,t_{_K})\right)$ may be defined accordingly.

\section{Block Space-Time Least-Squares Neural Network Method}

Our numerical results show that it is difficult to train the space-time LSNN method for problems with shock formation when the computational domain $\Omega$ is relatively large, even though the NN model we used is relatively small for approximating the solution of the underlying problem well. This is due to the nature of nonlinear hyperbolic conservation laws since information is transported from initial and inflow boundary to the rest of the domain along the flow direction. To overcome this training difficulty, we propose the block space-time LSNN method. 

\begin{figure}\label{domain_sketch}
    \centering
    \includegraphics[width=4.2in]{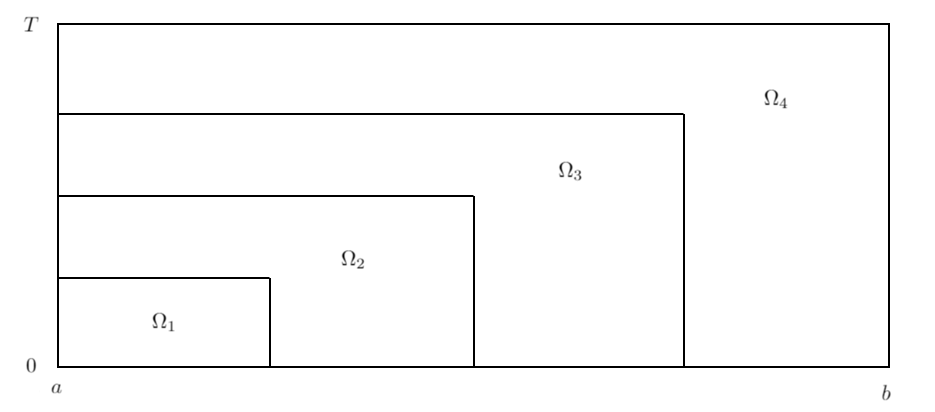}
    \caption{Sketch of domains $\{ \Omega_i\}_{i=1}^{m_0}$ for $m_0=4$}
    \label{fig:my_label}
\end{figure}

For clarity of exposition, let us consider one-dimensional problem defined on $\Omega= (a,b)\times (0,T)$. Without loss of generality, assume that $\tilde{\Gamma}_-=\{(a, t) |\, t\in (0, T)\}$ is the part of the boundary where the characteristic curves enter the domain $\Omega$. Hence, 
 \[
 \Gamma_-= \tilde{\Gamma}_- \cup \{(x, 0) |\, x\in (0, T)\}
\]
is the ``inflow'' boundary of $\Omega$. Let $m_0$ be a positive integer and let 
 \[
 \Omega_1=\left(a, a+\dfrac{b-a}{m_0}\right) \times \left(0, \dfrac{T}{m_0}\right)
 \quad\mbox{and }\,\, \Omega_i=\left(a, a+i\dfrac{(b-a)}{m_0}\right)\times \left(0, i\dfrac{T}{m_0}\right)\setminus \Omega_{i-1}
 \]
for $i=2, ..., m_0$. The sketch of domains $\{ \Omega_i\}_{i=1}^{m_0}$ is presented in Fig \ref{domain_sketch} and it is clear that $\{ \Omega_i\}_{i=1}^{m_0}$ forms a partition of the domain $\Omega$. 
Denote by $u_i=u|_{\Omega_i}$ the restriction of the solution $u$ of (\ref{pde}) on $\Omega_i$, then $u_i$ is the solution of the following problem:
\begin{equation} \label{pde2}
    \left\{\begin{array}{rccl}
    (u_i)_t+\nabla_{{\bx}} \!\cdot {\bff} (u_i) &= & 0, &\text{ in }\,\, {\Omega_i}\in \R^{2}, \\[2mm]
     u_i&=&
    {g},&\text{ on }\,\, \Gamma^i_-={\Gamma}_{-}\cap \partial \Omega_i ,
    \\[2mm]
    u_i&=&
    u_{i-1},&\text{ on }\,\, \Gamma_{i-1,i}=\partial\Omega_{i-1}\cap \partial \Omega_i
    \end{array}\right.
\end{equation}
for $i=1, ..., m_0$, where $\partial\Omega_0=\emptyset$. 

Define the least-squares functional for problem (\ref{pde2}) by
 \[
 \mathcal{L}^i\big(v;u_{i-1},g\big) = 
 \| v_t(\bx,t) + \nabla_{\bx} \!\cdot \bff (v)\|_{0,\Omega_i}^2 + \|v-g\|_{0, \Gamma^i_-}^2+ \|v-u_{i-1}\|_{0, \Gamma_{i-1,i}}^2.
 \]
Then the corresponding discrete least-squares functional
$\mathcal{L}^i_{_{\small {\cal T}}}\big(v(\cdot; {\small\btheta});u_{i-1},g\big)$ over the subdomain $\Omega_i$ may be defined in a similar fashion as in (\ref{L-NN-d}). Now, the block space-time LSNN method is to find 
${u}^i_{_{\small {\cal T}}}(\bzz,{\small\btheta}^*_i)\in \cM({\small\btheta},L)$ such that
 \begin{equation}\label{discrete_minimization_functional-block}
  \mathcal{L}^i_{_{\small {\cal T}}} \big({u}^i_{_{\small {\cal T}}}(\cdot,{\small\btheta}_i^*);u_{i-1},{g}\big) 
  = \min\limits_{v\in \cM({\scriptsize\btheta},L)} \mathcal{L}^i_{_{\small {\cal T}}}\big(v(\cdot;{\small\btheta});\,u_{i-1},{g}\big)
 = \min_{{\scriptsize \btheta}\in\R^{N}}\mathcal{L}^i_{_{\small {\cal T}}} \big(v(\cdot; {\small\btheta});u_{i-1},{g}\big)
\end{equation}
for $i=1,...,m_0$. 

\begin{remark}
The NN model $\cM({\scriptsize\btheta},L)$ is determined by the first subdomain and will be used for all subdomains. The trained parameter ${\small\btheta}^*_i$ from the $i^{th}$-subdomain is a good approximation to the parameters of the $(i+1)^{th}$-subdomain and, hence, may be used as an initial. This is because the solution in the current block is the evolution of the solution in the previous block.
\end{remark}

The block space-time LSNN method is based on a proper partition of the domain $\Omega$ depending on the ``inflow'' boundary of the domain. For example, in one dimension again, if 
\[
\Gamma_-=\{(x,t)\in [a,b]\times [0,T]|\, x=a, \,\,x=b, \mbox{ or }, t=0\},
\]
then the domain $\Omega$ may be partitioned by time blocks as 
\begin{equation}\label{1d_block_partition}
\Omega_i = (a,b)\times \big((i-1)T/m_0, iT/m_0\big)
\end{equation}
for $i=1,...,m_0$. Then the block space-time LSNN method may be defined accordingly.

\section{Implementation and Numerical Experiments}

In this section, we present numerical results for one dimensional benchmark test problems. Test problems include scalar nonlinear hyperbolic conservation law: (1) inviscid Burgers equation, i.e., ${\bff}(u) = \frac12 u^2$ (section 5.1-5.2) and (2) ${\bff}(u) = \frac14 u^4$ (section 5.3). Additionally, we analyze the effects of integration mesh and network structure in section 5.4; and compare the Roe and ENO schemes in section 5.5.

The domain $\Omega = (a,b)\times (0,T)$ is partitioned into time blocks as (\ref{1d_block_partition}) and $m_0$ is the number of blocks. Unless otherwise stated, the integration mesh $\mathcal{T}$ is obtained by uniformly partitioning all subdomains $\Omega_i$ into identical squares with the mesh size $h=0.01$ for $i=1,\cdots, m_0$. To preserve the conservation, the spatial mesh size $h_i$ and the temporal step size $\tau$ in both Roe (\ref{Roe}) and ENO (\ref{n-flux}) schemes are chosen to be the same as the integration mesh size, i.e., $h_i=\tau=h=0.01$. The block space-time LSNN method is implemented, and the minimization problem in (\ref{discrete_minimization_functional-block}) is numerically solved using the Adams version of gradient descent \cite{kingma2015} with a fixed or an adaptive learning rate.

As shown in \cite{Cai2021linear}, a three-layer NN is needed in order to approximate discontinuous solutions along non-straight line interfaces. Since there is no prior knowledge of geometric shape of the discontinuous interfaces for nonlinear problems, we use three-layer NN models for all test problems. Moreover, the same architecture of three-layer NN models is used for all blocks. As suggested in Remark 4.1, the parameters of the NN for the current block are initialized by the values of the parameters of the NN in the previous block. For the first block, the parameters of the second hidden layer are initialized randomly; and those of the first hidden layer are initialized using the strategy introduced in \cite{LiuCai1}. For convenience of readers, we briefly describe here. Let $\bomega_i\in \cS^{1}$ and $b_i\in \R$ be the weights and bias of the $i^{th}$ neuron of the first hidden layer of the first block NN model, respectively, where $\cS^{1}$ is the unit circle in $\R^2$. Initial of $\{(\bomega_i, b_i)\}_{i=1}^{n_1}$ is chosen so that the hyper-planes $\{\bomega_i\cdot(x,t)= b_i)\}_{i=1}^{n_1}$ form a uniform partition of the first block $\Omega_1$. In addition, without an effective training strategy, we observe from the experiment that adding a weight $\alpha$ to the $L^2$ loss of the initial condition is helpful for the training. Specifically, the following least-squares functional is used in the implementation
\begin{equation}\label{training_weight}
 \mathcal{L}^i\big(v;u_{i-1},g\big) = 
 \| v_t(\bx,t) + \nabla_{\bx} \!\cdot \bff (v)\|_{0,\Omega_i}^2 + \|v-g\|_{0, \Gamma^i_-}^2+ \alpha\|v-u_{i-1}\|_{0, \Gamma_{i-1,i}}^2 \quad \text{for } i=1,\cdots, m_0
\end{equation}
and the weight $\alpha$ is empirically determined.

Let $u_i$ be the solution of the problem in (\ref{pde2}) and $u_{i,_\cT}$ be the NN approximation. Using all quadrature points including those near the discontinuity, the relative error in the $L^2$ norm for each block is reported in Table \ref{riemann_shock_roe_table}-\ref{continuous_table}. Note that the points near discontinuity are usually excluded when reporting the error of existing traditional methods. The network structure is expressed as 2-$n_1$-$n_2$-1 for a three-layer network with $n_1$ and $n_2$ neurons in the respective first and second layers. The traces of the exact solution and the numerical approximation are depicted in Figures \ref{riemann_shock_roe_figure}-\ref{continuous_eno_figure} on a plane perpendicular to the space-time plane. Those traces exhibit the capability of the numerical approximation in resolving the shock/rarefaction. Since those planes are generally not perpendicular to the discontinuous interface, the errors shown in those traces are larger than the actual error.

\subsection{Riemann problem for the inviscid Burgers equation}
For the one dimensional inviscid Burgers equation, ${\bff}(u) =f(u)= \frac12 u^2$, we report numerical results for the corresponding Riemann problem where the initial condition with a single discontinuity is given by:
\begin{equation}\label{riemann_initial}
 u_0(x)
 =\left\{\begin{array}{rclll}
 u_{_L}, & \mbox{if }  x\leq 0,\\[2mm]
 u_{_R}, & \mbox{if } x >0.
 \end{array}\right.
\end{equation}
When $u_{_L}>u_{_R}$, the characteristic lines intersect and a shock forms immediately for $t > 0$. The weak solution is given by
\begin{equation}\label{riemann_shock_solution}
 u(x,t)
 =\left\{\begin{array}{rclll}
 u_{_L}, & \mbox{if }  x\leq st,\\[2mm]
u_{_R}, & \mbox{if } x >st,
 \end{array}\right.
\end{equation}
with the shock speed determined by the RH condition
\[s = \dfrac{f(u_{_L})- f(u_{_R})}{u_{_L}-u_{_R}}.
\]

When $u_{_L}<u_{_R}$, the range of influence of all points in $\R$ is a proper subset of $\R\times [0,\infty)$. This fact implies that the weak solution of the scalar hyperbolic conservation law is not unique. To ensure the underlying Cauchy problem having a unique solution over the whole domain $\R\times [0,\infty)$, the so-called vanishing viscosity weak solution is introduced (see, e.g., \cite{godlewski2013numerical,leveque1992numerical, thomas2013numerical}) and given by
\[
 u(x,t)
 =\left\{\begin{array}{rclll}
 u_{_L}, \quad  \mbox{if }  & x< u_{_L}t,\\[2mm]
 x/t, \quad \mbox{if }  & u_{_L}t\leq x\leq u_{_R}t,\\[2mm]
 u_{_R}, \quad \mbox{if } & x >u_{_R}t.
 \end{array}\right.
\]

\subsubsection{Shock formation}

The first test problem is corresponding to the case 
\[
u_{_L}=1> 0=u_{_R}
\]
with a computational domain $\Omega =(-1,2) \times (0,1)$. The inflow boundary is
\[
\Gamma_- = \Gamma_-^L\cup \Gamma_-^R\equiv \{(-1,t): t\in [0,1] \} \cup  \{(2,t): t\in [0,1]\}
\] 
with the boundary conditions: $g=u_{_L}$ on $\Gamma_-^L$ and $g=u_{_R}$ on $\Gamma_-^R$. The block space-time LSNN method is employed with $m_0=5$ blocks, a fixed learning rate 0.003, and 30000 iterations for each block.

\begin{figure}[htbp]
  \centering 
  \subfigure[Network approximation $u_{_\cT}$ \newline on $\Omega$]{ 
    \includegraphics[width=1.8in]{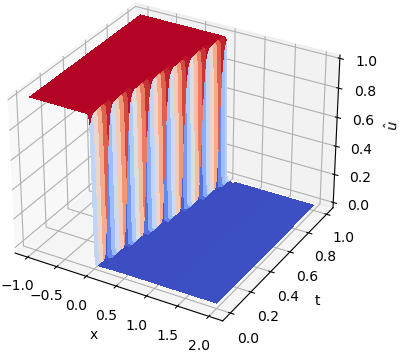}}
    \subfigure[Traces of exact and numerical \newline solutions $u_{1,_\cT}$ on the plane $t=0.2$]{ 
    \includegraphics[width=1.88in]{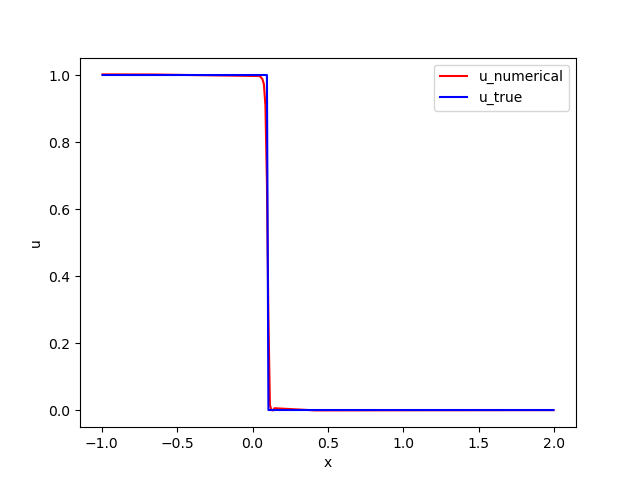}} 
    \subfigure[Traces of exact and numerical \newline solutions $u_{2,_\cT}$ on the plane $t=0.4$]{ 
    \includegraphics[width=1.88in]{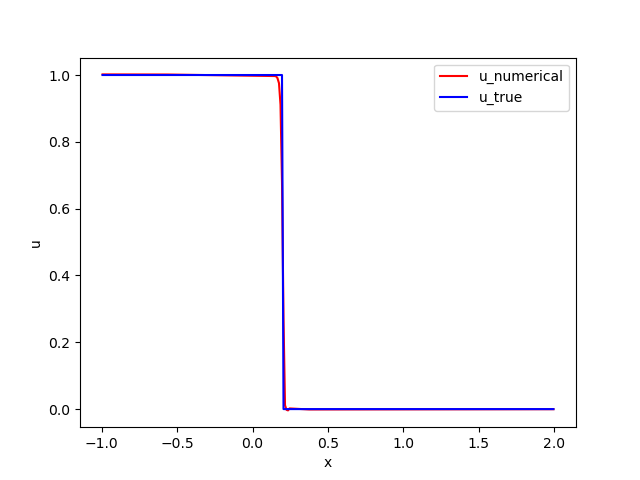}} 
\\
  \subfigure[Traces of exact and numerical \newline solutions $u_{3,_\cT}$ on the plane $t=0.6$]{ 
    \includegraphics[width=1.88in]{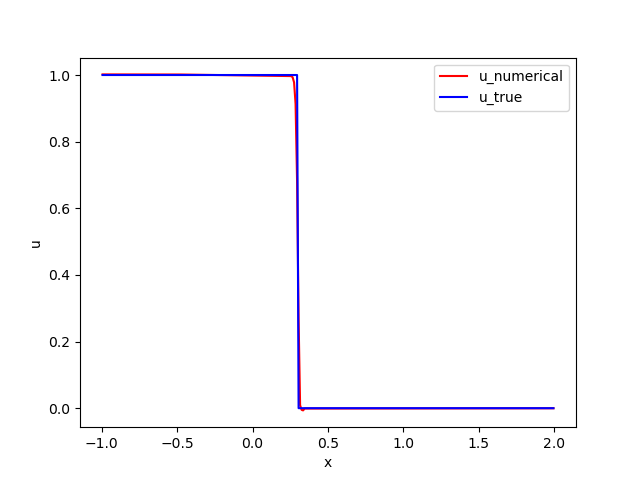}}
      \subfigure[Traces of exact and numerical \newline solutions $u_{4,_\cT}$ on the plane $t=0.8$]{ 
    \includegraphics[width=1.88in]{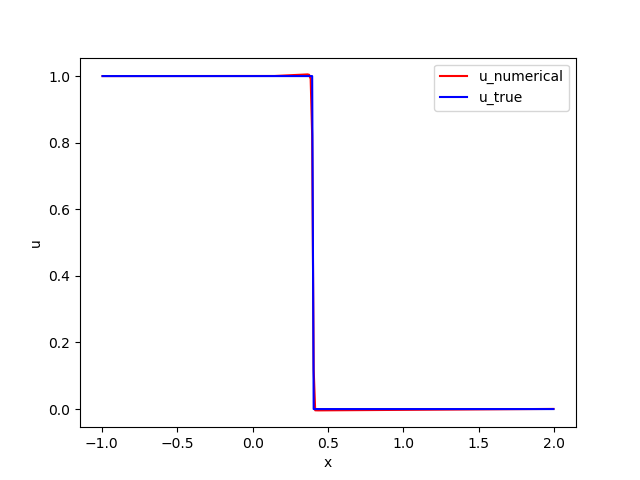}}
      \subfigure[Traces of exact and numerical \newline solutions $u_{5,_\cT}$ on the plane $t=1.0$]{ 
    \includegraphics[width=1.88in]{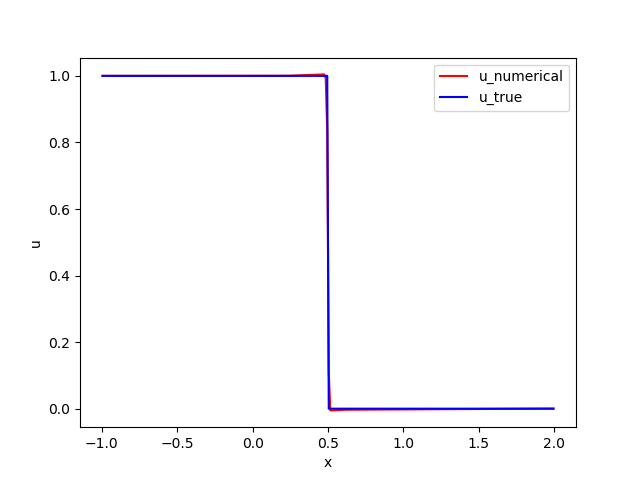}}
  \caption{Approximation results of Riemann problem (shock) for Burgers' equation using Roe flux} 
  \label{riemann_shock_roe_figure}
\end{figure}

The set of experiments is done by using the numerical fluxes of Roe (\ref{Roe}) and the second order ENO (\ref{n-flux}). By choosing $\alpha =20$ in (\ref{training_weight}), numerical results of both schemes are reported in Table \ref{riemann_shock_roe_table}. Since the results of the Roe flux are similar, Figure \ref{riemann_shock_roe_figure} (b)-(f) only depict the traces of the exact and numerical solution generated by the ENO flux on the planes $t=iT/m_0$ for $i=1,\cdots, m_0$. Clearly, the block space-time LSNN method with a conservative scheme is able to resolve the shock and accurately approximate the solution without the Gibbs phenomena. For this simple Riemann problem, Roe and ENO schemes produce similar results. Given the additional evaluations involving in the ENO scheme, we do not observe many advantages of using higher-order scheme for flux reconstruction despite the fact that ENO performs slightly better in the $L^2$ relative errors.

 \begin{table}[htbp]
\centering
\caption{Relative errors of Riemann problem (shock) for Burgers' equation using Roe and ENO fluxes}
\vspace{5pt}
\begin{tabular}{|l |c | c | c | c |  c| c|c|}
\hline
 Network structure &Block &
\multicolumn{1}{|p{4cm}|}{\centering Roe flux \\ $\frac{\|u_i-u_{i,_\cT}\|_0}{\|u_i\|_0}$ } &
\multicolumn{1}{|p{4cm}|}{\centering ENO flux \\ $\frac{\|u_i-u_{i,_\cT}\|_0}{\|u_i\|_0}$ } \\ \hline
2-10-10-1 & $\Omega_1$  &0.049553  & 0.030901 \\ \hline
2-10-10-1 & $\Omega_2$  &0.046321  & 0.030178\\\hline
2-10-10-1 & $\Omega_3$  &0.044123   &0.028984 \\\hline
2-10-10-1 & $\Omega_4$  &0.042621 &0.028791  \\ \hline
2-10-10-1 & $\Omega_5$  & 0.041182 &0.032253 \\ \hline
\end{tabular}
\centering
\label{riemann_shock_roe_table}
\end{table}

\subsubsection{Rarefaction waves}

The second test problem is corresponding to the case 
\[
u_{_L}=0 < 1=u_{_R}
\]
with a computational domain $\Omega =(-1,2) \times (0,0.2)$. The inflow boundary is
\[
\Gamma_- = \Gamma_-^L\cup \Gamma_-^R\equiv \{(-1,t): t\in [0,0.2] \} \cup  \{(2,t): t\in [0,0.2]\}
\] 
with the boundary conditions: $g=u_{_L}$ on $\Gamma_-^L$ and $g=u_{_R}$ on $\Gamma_-^R$. The block space-time LSNN method is employed with $m_0=2$ blocks, a fixed learning rate 0.003, and 20000 iterations for each block.

\begin{figure}[htbp]
  \centering 
  \subfigure[Exact solution $u$ on $\Omega$]{
    \includegraphics[width=1.65in]{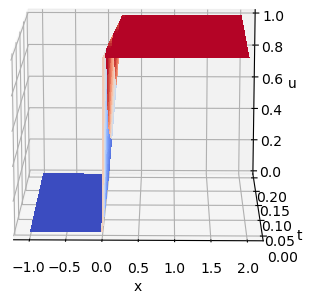}}
      \hspace{0.06in} 
      \subfigure[Traces of exact and numerical \newline solutions $u_{1,_\cT}$ on the plane $t=0.1$]{
    \includegraphics[width=1.9in]{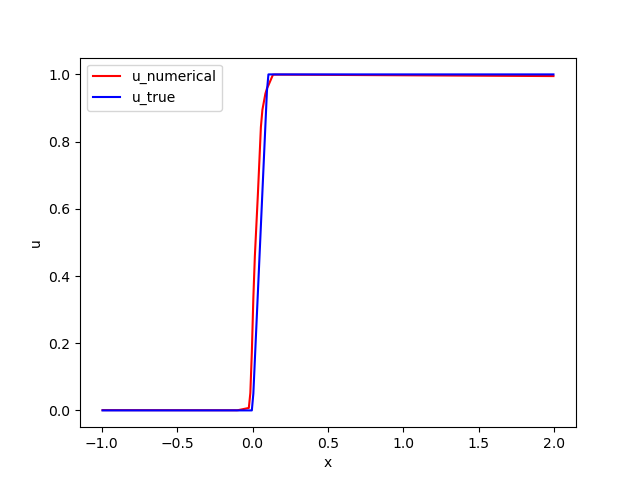}}
    \hspace{0.06in} 
      \subfigure[Traces of exact and numerical \newline solutions $u_{2,_\cT}$ on the plane $t=0.2$]{
    \includegraphics[width=1.9in]{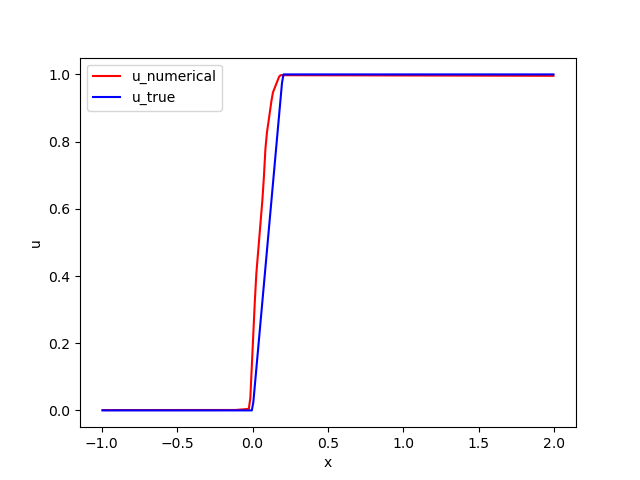}}
  \caption{Approximation results of Riemann problem (rarefaction) for Burgers' equation using Roe flux} 
  \label{riemann_rare_roe_figure}
\end{figure}

Numerical results of a 2-10-10-1 network using the Roe flux (\ref{Roe}) are reported in Table \ref{riemann_rare_roe_table}. The traces of the exact and numerical solutions in Fig. \ref{riemann_rare_roe_figure} indicate that Roe's scheme fails to resolve the rarefaction. For the traditional mesh-based method, it is well-known that Roe's scheme is not able to compute the physical solution of the rarefaction wave. This is because the scheme approximates the numerical flux depending the sign of the speed $a_i$ (\ref{roe_a}) at midpoints. If the sign differs on two sides and $u(x,t)$ travels slower on the left, Roe's scheme may not be able to capture such behavior. From this perspective, we observe certain limitations of using such conservative scheme and that the discretization scheme is very important for the block space-time LSNN method.

\begin{remark}
In \em \cite{harten1997high, leveque1992numerical}, the authors proposed the ``entropy fix'' traditional mesh-based approach to address such issue. In a forthcoming paper, we will propose a discretization scheme for the LSNN method which is capable of resolving the rarefaction.
\end{remark}

\begin{table}[htbp]
\centering
\caption{Relative errors of Riemann problem (rarefaction) for Burgers' equation using Roe flux}
\vspace{5pt}
\begin{tabular}{|l|l|l|}
\hline
Network structure &Time block & $\frac{\|u_i-u_{i,_\cT}\|_0}{\|u_i\|_0}$ \\ \hline
2-10-10-1 & $\Omega_1 $   & 0.047435\\ \hline
2-10-10-1 & $\Omega_2 $  & 0.074521\\ \hline
\end{tabular}
\centering
\label{riemann_rare_roe_table}
\end{table}

\subsection{Inviscid Burgers equation with smooth initial condition}
The third problem is again the Burgers equation defined on the computational domain $\Omega = (0,2)\times (0,0.4)$ with the inflow boundary 
\[
\Gamma_- = \Gamma_-^L\cup \Gamma_-^R\equiv \{(0,t): t\in [0,0.4] \} \cup  \{(2,t): t\in [0,0.4]\}
\]  
and a sinusoidal initial condition
\[u_0(x) = 0.5 + \sin(\pi x).
\]
The shock forms at $t = 1/\pi \approx 0.318$. Since the exact solution of the test problem is defined implicitly, in order to measure the quality of the NN approximation, we generate a benchmark reference solution $\hat{u}$ using the traditional mesh-based approach. Specifically, the third order accurate WENO scheme \cite{shu1998essentially} is employed for the spatial discretization with a fine grid ($\Delta x = 0.001$ and $\Delta t =0.0002$) on the computational domain $\Omega$; and the fourth order Runge-Kutta method is used for the temporal discretization \cite{wang2019learning}. The block space-time LSNN method is implemented with $m_0=8$ blocks and an adaptive learning rate which starts at 0.005 and decays by half for every 25000 iterations. The learning rate decay strategy is employed with 50000 iterations on each time block.

\begin{figure}[htbp]
  \centering 
    \subfigure[Traces of reference and numerical solutions $u_{1,_\cT}$ on the plane $t=0.05$]{ 
    \includegraphics[width=1.7in]{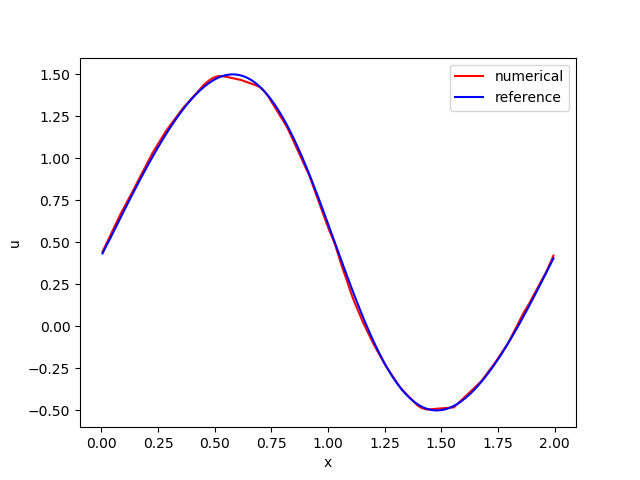}} 
  \hspace{0.2in} 
  \subfigure[Traces of reference and numerical solutions $u_{2,_\cT}$ on the plane $t=0.1$]{ 
    \includegraphics[width=1.7in]{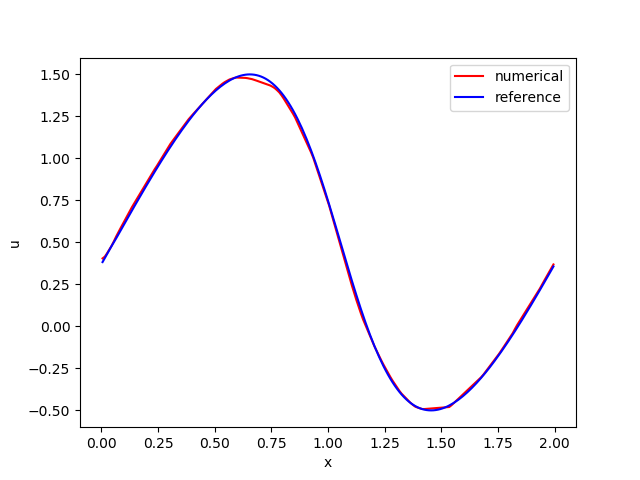}} 
  \hspace{0.2in} 
  \subfigure[Traces of reference and numerical solutions $u_{3,_\cT}$ on the plane $t=0.15$]{ 
    \includegraphics[width=1.7in]{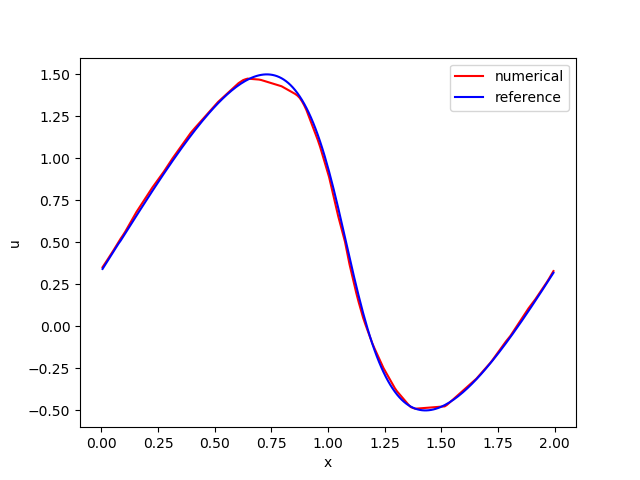}}
\\
      \subfigure[Traces of reference and numerical solutions $u_{4,_\cT}$ on the plane $t=0.2$]{ 
    \includegraphics[width=1.7in]{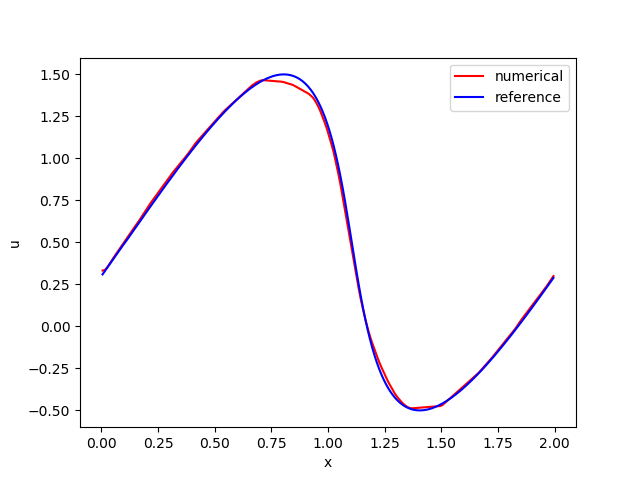}}
     \hspace{0.2in} 
\subfigure[Traces of reference and numerical solutions $u_{5,_\cT}$ on the plane $t=0.25$]{ 
    \includegraphics[width=1.7in]{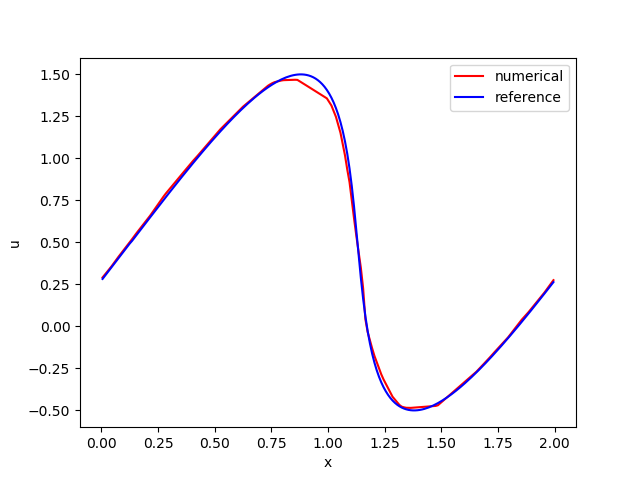}} 
    \hspace{0.2in} 
\subfigure[Traces of reference and numerical solutions $u_{6,_\cT}$ on the plane $t=0.3$]{ 
    \includegraphics[width=1.7in]{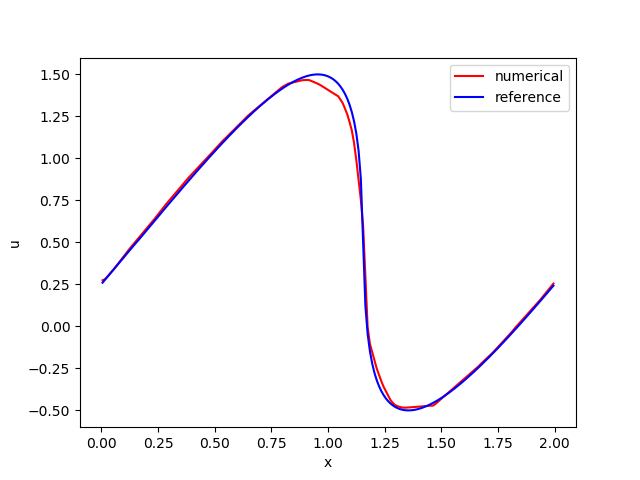}}
\\
     \subfigure[Traces of reference and numerical solutions $u_{7,_\cT}$ on the plane $t=0.35$]{ 
    \includegraphics[width=1.7in]{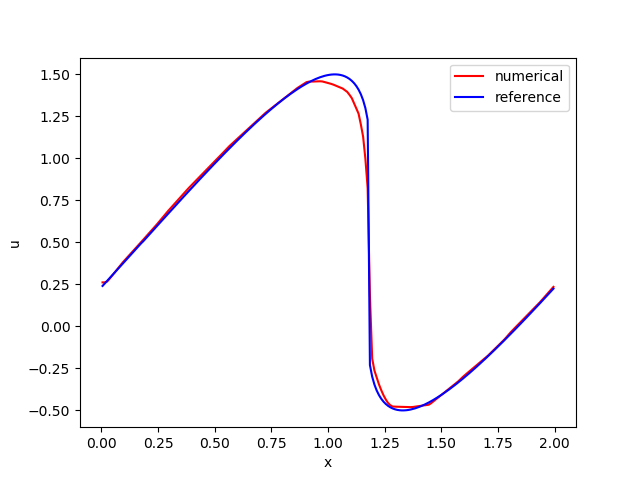}}
 \hspace{0.2in} 
\subfigure[Traces of reference and numerical solutions $u_{8,_\cT}$ on the plane $t=0.4$]{ 
    \includegraphics[width=1.7in]{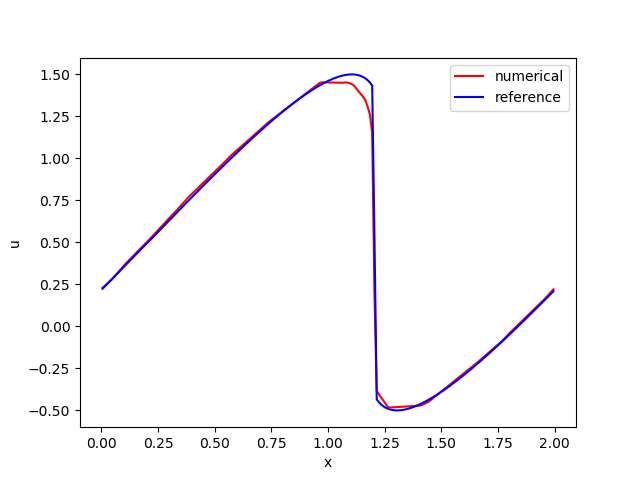}}
  \caption{Approximation results of Burgers' equation with a sinusoidal initial using ENO flux} 
  \label{sin_eno_figure}
\end{figure}

Since the initial condition of the test problem is a smooth function, it is expected that a network with additional neurons is needed for approximation. Choosing $\alpha=5$ in (\ref{training_weight}), the numerical results of a 2-30-30-1 network using the ENO flux in (\ref{n-flux}) are reported in Table \ref{sin_eno_table}. Figure \ref{sin_eno_figure} depicts the traces of the reference solution and numerical approximation on the plane $t=iT/m_0$ for $i=1,\cdots, m_0$. We observe some error accumulation when block evolves, and the block space-time LSNN method can resolve the shock. It is noticeable in Fig. \ref{sin_eno_figure} that approximation near the local maximum is poor. Possibly this is due to inaccuracy of the second order ENO scheme for complicated initial data. A new and accurate discretization scheme for the LSNN method will be reported in a forthcoming paper \cite{Cai2021FV}.

\begin{table}[htbp]
\centering
\caption{Relative errors of Burgers' equation with a sinusoidal initial condition using ENO flux}
\vspace{5pt}
\begin{tabular}{|l|l|l|}
\hline
Network structure &Time block & $\frac{\|\hat{u}_i-u_{i,_\cT}\|_0}{\|\hat{u}_i\|_0}$ \\ \hline
2-30-30-1 & $\Omega_1$ &0.010461 \\ \hline
2-30-30-1 & $\Omega_2$ &0.012517 \\ \hline
2-30-30-1 & $\Omega_3$ & 0.019772\\ \hline
2-30-30-1 & $\Omega_4$ & 0.022574\\ \hline
2-30-30-1 & $\Omega_5$ & 0.029011\\ \hline
2-30-30-1 & $\Omega_6$ & 0.038852\\ \hline
2-30-30-1 & $\Omega_7$ & 0.075888\\ \hline
2-30-30-1 & $\Omega_8$ &0.078581 \\ \hline
\end{tabular}
\centering
\label{sin_eno_table}
\end{table}

\subsection[alternative title goes here]{Riemann problem with $f(u)=\frac14 u^4$}
The fourth numerical experiment is the Riemann shock problem with a convex flux ${\bff}(u) = f(u)=\frac14 u^4$. We choose the initial condition
 \[
 u_{_L}=1 >0=u_{_R}
 \]
 in (\ref{riemann_initial}), then the weak solution is given by (\ref{riemann_shock_solution}) with the speed $s=1/4$. The computational domain of problem is given by $\Sigma = (-1,1)\times (0,1)$ and the inflow boundary is
\[
\Gamma_- = \Gamma_-^L\cup \Gamma_-^R\equiv \{(-1,t): t\in [0,1] \} \cup  \{(1,t): t\in [0,1]\}
\] 
with the boundary conditions: $g=1$ on $\Gamma_-^L$ and $g=0$ on $\Gamma_-^R$. 

\begin{figure}[htbp]
  \centering 
  \subfigure[Network approximation $u_{_\cT}$ on $\Omega$]{ 
    \includegraphics[width=1.98in]{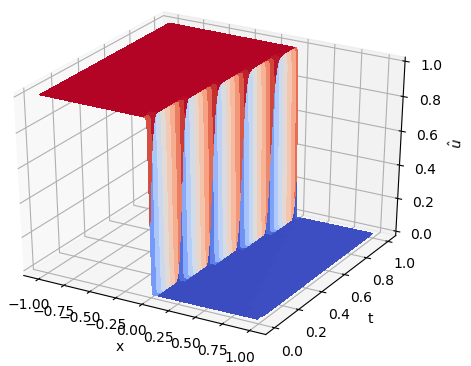}}
\hspace{0.3in}
    \subfigure[Traces of exact and numerical solutions $u_{5,_\cT}$ on the plane $t=1.0$]{ 
    \includegraphics[width=2in]{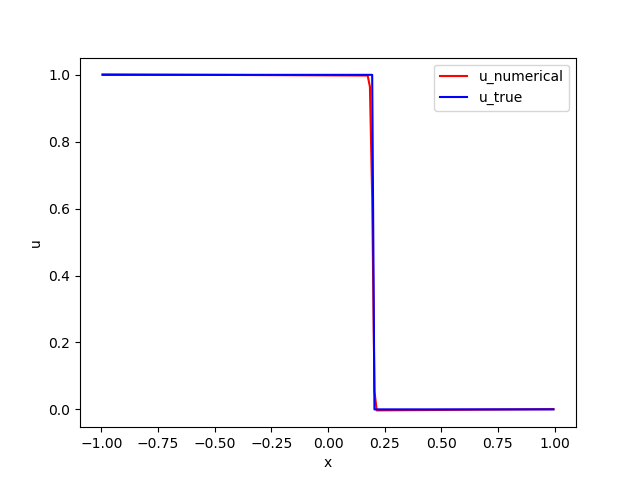}} 
  \caption{Approximation results of Riemann problem (shock) with $f(u) = \frac14 u^4$ using Roe flux} 
  \label{convex_shock_figure}
\end{figure}

Empirically, we choose $\alpha=20$ in (\ref{training_weight}) in the implementation. Employing the block space-time LSNN method with $m_0=5$ blocks, a fixed learning rate 0.003 and 30000 iterations for each block, we report the numerical results of a 2-10-10-1 network in Table \ref{convex_shock_table} and Fig. \ref{convex_shock_figure}. Obviously, the discontinuous interface can be accurately captured as the NN approximation is almost overlapped with the exact solution,  which suggests that the LSNN method is not only capable of solving the Burgers equation but also the problem with a general convex flux.
\begin{table}[htbp]
\centering
\caption{Relative errors of Riemann problem (shock) with $f(u) = \frac14 u^4$ using Roe flux}
\vspace{5pt}
\begin{tabular}{|l|l|l|}
\hline
Network structure &Block & $\frac{\|u_i-u_{i,_\cT}\|_0}{\|u_i\|_0}$ \\ \hline
2-10-10-1 & $\Omega_1$  &0.035034  \\ \hline
2-10-10-1 & $\Omega_2$  &0.036645  \\ \hline
2-10-10-1 & $\Omega_3$  & 0.036798 \\ \hline
2-10-10-1 & $\Omega_4$  &0.037217  \\ \hline
2-10-10-1 & $\Omega_5$  & 0.037451 \\ \hline
\end{tabular}
\centering
\label{convex_shock_table}
\end{table}

\subsection{Effects of integration mesh and network structure}
The goal of this section is to analyze the effects of integration mesh and network structure for the LSNN method. Specifically, we use the inviscid Burgers equation defined on the computational domain $ \Omega  =(-1,2) \times (0,0.6)$ with a continuous piece-wise linear initial condition
\[
 u_0(x)
 =\left\{\begin{array}{rclll}
 1, &\quad\mbox{if }&  x< 0,\\[2mm]
 1-2x, &\quad\mbox{if }&  0\leq x \leq 1/2,\\[2mm]
 0, &\quad\mbox{if }& x >1/2.
 \end{array}\right.
\]
The inflow boundary is
\[
\Gamma_- = \Gamma_-^L\cup \Gamma_-^R\equiv \{(-1,t): t\in [0,0.6] \} \cup  \{(2,t): t\in [0,0.6]\},
\] 
with the boundary conditions: $g=1$ on $\Gamma_-^L$ and $g=0$ on $\Gamma_-^R$. Even though the initial value $u_0$ is continuous, the shock will appear at some point since $u(x,t)$ travels faster on the left-hand side than on the right-hand side. Specifically, when $t<1/2$, the solution is continuous and it is determined by the characteristic lines as well as the initial conditions:
\[
 u(x,t)
 =\left\{\begin{array}{rclll}
 1, & \mbox{if }&  x< t<1/2,\\[2mm]
 \dfrac{1-2x}{1-2t}, & \mbox{if }& t\leq x \leq 1/2,\\[2mm]
 0, & \mbox{if }& x >1/2.
 \end{array}\right.
\]
 When $t>1/2$, the shock forms and the desired weak solution satisfying RH condition is given by
\[
 u(x,t)
 =\left\{\begin{array}{rclll}
 1, & \mbox{if }&  x<(2t+1)/4,\\[2mm]
 0, & \mbox{if }& x \geq (2t+1)/4.
 \end{array}\right.
\]

The block space-time LSNN is implemented with $m_0=3$ blocks, a fixed learning rate 0.003 and $\alpha=10$ in the training. In order to explore the network approximation power, we do not constrain the number of iterations on each block. The stopping criteria for the gradient descent solver is set as follows: the solver stops when the loss function (\ref{training_weight}) decreases within $0.1\%$ in the last 2000 iterations.

The first set of experiments is to observe the impact of integration mesh size when a fixed network structure 2-10-10-1 is used. Starting with the same initialization for both layers, the relative $L^2$ errors of the ENO scheme on uniform meshes with different mesh size are reported in Table \ref{quadrature_convergence}. The results display that the error decreases as the integration mesh size reduces, but the decreasing rate is small when the mesh size is sufficiently fine. This indicates that $h=0.01$ is fine enough to accurately evaluate the discrete LS functional in (\ref{discrete_minimization_functional-block}) for the given 2-10-10-1 network.


\begin{table}[htbp]
\centering
\caption{Relative $L^2$ errors of the problem with a piece-wise linear initial using different integration mesh sizes}
\vspace{5pt}
\begin{tabular}{|c|l|l|l|l|l|}
\hline
\multirow{2}{*}{Time block} & \multicolumn{5}{c|}{Integration mesh size} \\ \cline{2-6} 
                            & $h=0.1$  & $h=0.05$ & $h=0.025$  & $h=0.01$  & $h=0.005$ \\ \hline
\multicolumn{1}{|l|}{ $\Omega_1$} &0.056712 &0.029366 &0.017376 &0.013946& 0.013121 \\ \hline
\multicolumn{1}{|l|}{ $\Omega_2$} &0.079633 &0.064069 &0.048971 &0.037891& 0.036501 \\ \hline
\multicolumn{1}{|l|}{ $\Omega_3$} &0.120129 &0.115123 &0.102844 &0.084728& 0.081209  \\ \hline
\end{tabular}
\label{quadrature_convergence}
\end{table}

The second set of numerical tests is to investigate the effect of network structure. To accurately evaluate the functional, we use a fine integration mesh with size $h=0.005$. The approximation errors generated by four different network structures using the ENO scheme are presented in Table \ref{network_convergence}. As expected, the approximation error decreases as the number of neurons in the network increases.

\begin{table}[htbp]
\centering
\caption{Relative $L^2$ errors of the problem with a piece-wise linear initial using different network structures}
\vspace{5pt}
\begin{tabular}{|c|l|l|l|l|}
\hline
\multirow{2}{*}{Time block} & \multicolumn{4}{c|}{Network structure} \\ \cline{2-5} 
& 2-4-4-1  & 2-7-7-1  & 2-10-10-1  & 2-13-13-1 \\ \hline
\multicolumn{1}{|l|}{ $\Omega_1$}&0.021699 &0.017004 &0.013121 &0.009895\\ \hline
\multicolumn{1}{|l|}{ $\Omega_2$}&0.047953 &0.040521 &0.036501 &0.028159\\ \hline
\multicolumn{1}{|l|}{ $\Omega_3$}&0.095649 &0.088201 &0.081209  &0.072244\\ \hline
\end{tabular}
\label{network_convergence}
\end{table}

Tables \ref{quadrature_convergence} and \ref{network_convergence} indicate that the approximation error grows significantly as the block evolves. To address this issue, a simple way is to decrease the size of each block, or equivalently increase the number of blocks $m_0$. Specifically, increasing $m_0$ from 3 to 6, Table \ref{block_number} reports numerical results of the ENO scheme using $h=0.01$, which is better than the results in the corresponding column in Table \ref{quadrature_convergence}. Due to the increasing computational cost when using a smaller block, in practice we need to balance the cost and the accuracy when choosing the size of the block in the block space-time LSNN method.

\begin{table}[htbp]
\centering
\caption{Relative $L^2$ errors of the problem with a piece-wise linear initial using ENO flux with $m_0=6$}
\vspace{5pt}
\begin{tabular}{|l|l|l|}
\hline
Network structure &Block & $\frac{\|u_i-u_{i,_\cT}\|_0}{\|u_i\|_0}$ \\ \hline
 2-10-10-1 & $\Omega_1$  &0.008958  \\ \hline
 2-10-10-1 & $\Omega_2$  &0.014842  \\ \hline
 2-10-10-1 & $\Omega_3$  &0.019282  \\ \hline
 2-10-10-1 & $\Omega_4$  &0.025562  \\ \hline
 2-10-10-1 & $\Omega_5$  &0.047431  \\ \hline
 2-10-10-1 & $\Omega_6$  &0.052606  \\ \hline
\end{tabular}
\centering
\label{block_number}
\end{table}

\subsection{Comparison of the Roe and ENO schemes}
\begin{figure}[htbp]
  \centering 
    \subfigure[Traces of exact and numerical solutions $u_{2,_\cT}$ on the plane $t=0.2$]{
    \includegraphics[width=1.75in]{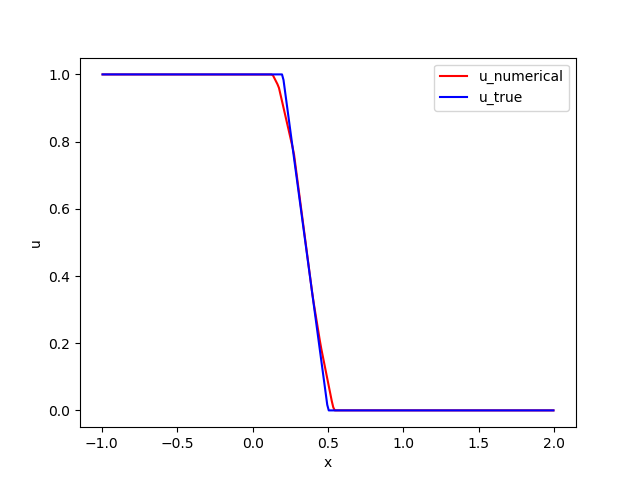}} 
  \hspace{0.1in} 
  \subfigure[Traces of exact and numerical solutions $u_{4,_\cT}$ on the plane $t=0.4$]{ 
    \includegraphics[width=1.75in]{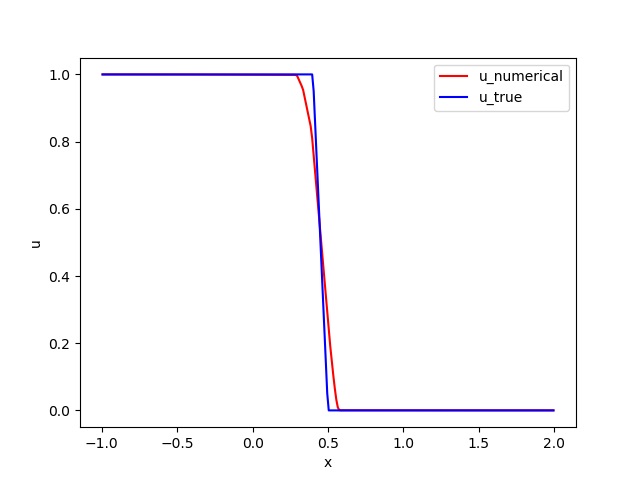}} 
  \hspace{0.1in} 
  \subfigure[Traces of exact and numerical solutions $u_{6,_\cT}$ on the plane $t=0.6$]{
    \includegraphics[width=1.75in]{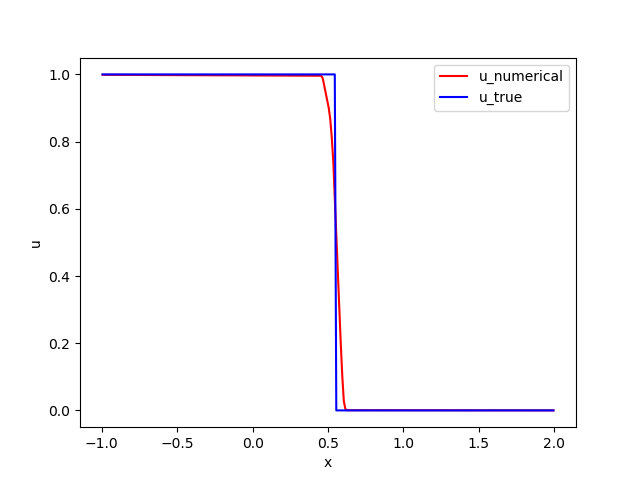}}
  \caption{Approximation results of Burgers' equation with a piece-wise linear initial using Roe flux} 
  \label{continuous_roe_figure}
\end{figure}

\begin{figure}[htbp]
  \centering 
    \subfigure[Traces of exact and numerical solutions $u_{2,_\cT}$ on the plane $t=0.2$]{ 
    \includegraphics[width=1.75in]{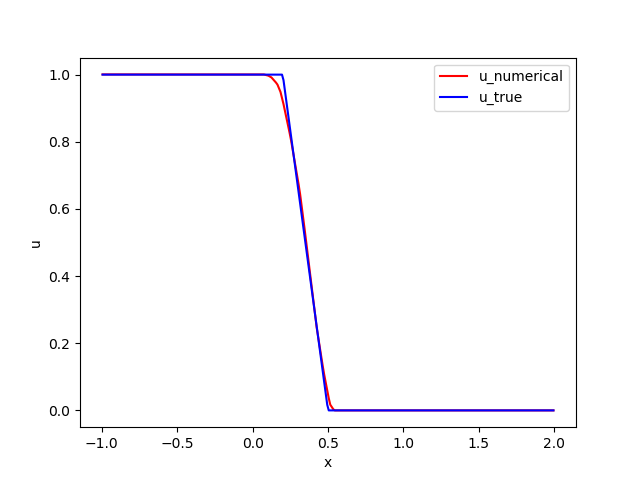}} 
  \hspace{0.1in} 
  \subfigure[Traces of exact and numerical solutions $u_{4,_\cT}$ on the plane $t=0.4$]{
    \includegraphics[width=1.75in]{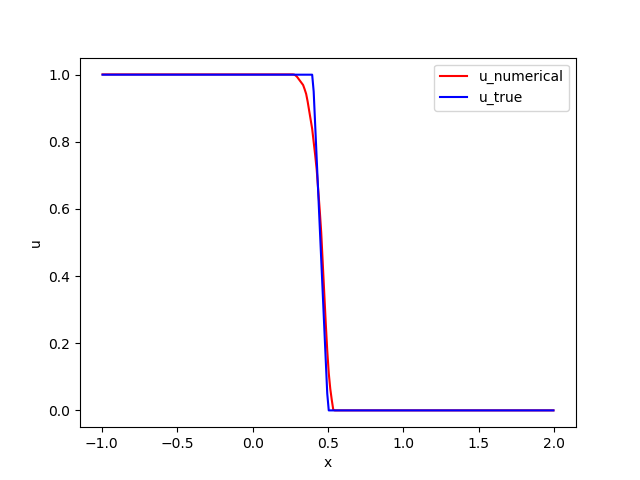}} 
  \hspace{0.1in} 
  \subfigure[Traces of exact and numerical solutions $u_{6,_\cT}$ on the plane $t=0.6$]{
    \includegraphics[width=1.75in]{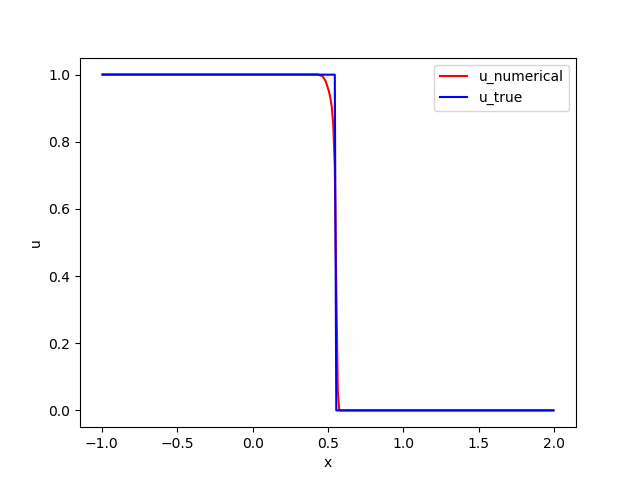}}
  \caption{Approximation results of Burgers' equation with a piece-wise linear initial using ENO flux} 
  \label{continuous_eno_figure}
\end{figure}

This choosing set of numerical experiments aims to compare numerical performances of the LSNN method using the Roe (\ref{Roe}) and ENO (\ref{n-flux}) schemes. The test problem in section 5.4 is used, and the block space-time LSNN method is employed with $m_0=6$ blocks, the 2-10-10-1 network structure, and with a uniform integration mesh of the size $h=0.01$. We note that the ENO performs better than the Roe in the relative $L^2$ norm (see Table \ref{continuous_table}) and near the discontinuous interface (see Fig. \ref{continuous_roe_figure} and \ref{continuous_eno_figure}). Note that the former is more expansive than the latter in the computational cost due to additional testing.

\begin{table}[htbp]
\centering
\caption{Relative errors of Burgers' equation with a piece-wise linear initial condition}
\vspace{5pt}
\begin{tabular}{|l |c | c | c | c |  c| c|c|c|c|c|}
\hline
 Network structure &Block &
\multicolumn{1}{|p{4cm}|}{\centering Roe flux \\ $\frac{\|u_i-u_{i,_\cT}\|_0}{\|u_i\|_0}$ } &
\multicolumn{1}{|p{4cm}|}{\centering ENO flux \\ $\frac{\|u_i-u_{i,_\cT}\|_0}{\|u_i\|_0}$ }
\\ \hline
2-10-10-1 & $\Omega_1$  &0.008333  &0.008958  \\ \hline
2-10-10-1 & $\Omega_2$  &0.017503  &0.014842 \\\hline
2-10-10-1 & $\Omega_3$  &0.025342  &0.019282 \\\hline
2-10-10-1 & $\Omega_4$  &0.036211  &0.025562  \\ \hline
2-10-10-1 & $\Omega_5$  &0.061251  &0.047431 \\\hline
2-10-10-1 & $\Omega_6$  &0.068512  &0.052606 \\\hline
\end{tabular}
\centering
\label{continuous_table}
\end{table}

\section{Discussions and Conclusions}

The block space-time LSNN method is proposed for solving scalar nonlinear hyperbolic conservation laws. The least-squares formulation is a direct application of the least-squares principle to the underlying problem: the equation, the inflow boundary condition, and the initial condition. The block space-time version of the LSNN method is introduced to compensate with some uncertainty of the not well-understood nonlinear optimization procedure. 

How to approximate the differential operator in the least-squares functional is critical for the success of the space-time LSNN method. As mentioned in the introduction, existing NN-based methods are not applicable to the inviscid Burgers equation whose solution is discontinuous. Employing the Roe and the second order ENO schemes, we show numerically that the resulting LSNN method is capable of resolving the shock without smearing and oscillations. Moreover, the LSNN method has much less degrees of freedom (DoF) than traditional mesh-based methods.

Despite the great potential demonstrated in this paper, the current version of the LSNN method needs to be improved in both accuracy and efficiency in order to grow into a viable numerical method. First, the method is inaccurate for complicated initial condition; moreover, it has limitations for problems with rarefaction waves and with non-convex spatial fluxes. To overcome these deficiency, we recently study a new version of the LSNN method that uses a novel and accurate finite volume method developed in \cite{Cai2021FV}. Second, the computational cost of the current version of the LSNN method is expensive because of the resulting non-convex optimization. 
To reduce the cost, we will explore physical meanings of non-linear parameters of NNs (see, e.g., \cite{LiuCai1}) to develop good initial strategies and then employ fast iterative solvers such as methods of BFGS type (the Broyden–Fletcher–Goldfarb–Shanno algorithm).



\bigskip
\bibliographystyle{ieee}
\bibliography{main.bbl}

\begin{thebibliography}{10}\itemsep=-1pt

\bibitem{PNAS2019}
Y.~Bar-Sinai, S.~Hoyer, J.~Hickey, and M.~P. Brenner.
\newblock Learning data-driven discretizations for partial differential
  equations.
\newblock {\em Proceedings of the National Academy of Science of USA}, 116
  (31):15344--15349, 2019.

\bibitem{bochev2001improved}
P.~Bochev and J.~Choi.
\newblock Improved least-squares error estimates for scalar hyperbolic
  problems.
\newblock {\em Comput. Methods Appl. Math.}, 1(2):115--124, 2001.

\bibitem{brezzi2004discontinuous}
F.~Brezzi, L.~D. Marini, and E.~S{\"u}li.
\newblock Discontinuous galerkin methods for first-order hyperbolic problems.
\newblock {\em Math. Models Methods Appl. Sci.}, 14(12):1893--1903, 2004.

\bibitem{burman2009posteriori}
E.~Burman.
\newblock A posteriori error estimation for interior penalty finite element
  approximations of the advection-reaction equation.
\newblock {\em SIAM J. Numer. Anal.}, 47(5):3584--3607, 2009.

\bibitem{Cai2021FV}
Z.~Cai, J.~Chen, and M.~Liu.
\newblock Finite volume least-squares neural network ({FV-LSNN}) method for
  scalar nonlinear hyperbolic conservation laws.
\newblock {\em arXiv preprint arXiv:2110.10895}, 2021.

\bibitem{Cai2021linear}
Z.~Cai, J.~Chen, and M.~Liu.
\newblock Least-squares {ReLU} neural network {(LSNN)} method for linear
  advection-reaction equation.
\newblock {\em J. Comput. Phys.}, 443 (2021) 110514.

\bibitem{cai2020deep}
Z.~Cai, J.~Chen, M.~Liu, and X.~Liu.
\newblock Deep least-squares methods: An unsupervised learning-based numerical
  method for solving elliptic {PDE}s.
\newblock {\em J. Comput. Phys.}, 420 (2020) 109707.

\bibitem{dahmen2012adaptive}
W.~Dahmen, C.~Huang, C.~Schwab, and G.~Welper.
\newblock Adaptive petrov--galerkin methods for first order transport
  equations.
\newblock {\em SIAM J. Numer. Anal.}, 50(5):2420--2445, 2012.

\bibitem{de2004least}
H.~De~Sterck, T.~A. Manteuffel, S.~F. McCormick, and L.~Olson.
\newblock Least-squares finite element methods and algebraic multigrid solvers
  for linear hyperbolic {PDE}s.
\newblock {\em SIAM J. Sci. Comput.}, 26(1):31--54, 2004.

\bibitem{de2005numerical}
H.~De~Sterck, T.~A. Manteuffel, S.~F. McCormick, and L.~Olson.
\newblock Numerical conservation properties of {H}(div)-conforming
  least-squares finite element methods for the burgers equation.
\newblock {\em SIAM J. Sci. Comput.}, 26(5):1573--1597, 2005.

\bibitem{demkowicz2010class}
L.~Demkowicz and J.~Gopalakrishnan.
\newblock A class of discontinuous petrov--galerkin methods. part {I}: The
  transport equation.
\newblock {\em Comput. Methods Appl. Mech. Eng.}, 199(23-24):1558--1572, 2010.

\bibitem{Weinan18}
W.~E and B.~Yu.
\newblock The deep ritz method: A deep learning-based numerical algorithm for
  solving variational problems.
\newblock {\em Communications in Mathematics and Statistics}, 6(1):1--12, 3
  2018.

\bibitem{godlewski2013numerical}
E.~Godlewski and P.-A. Raviart.
\newblock {\em Numerical Approximation of Hyperbolic Systems of Conservation
  Laws}, volume 118.
\newblock Springer Science \& Business Media, 2013.

\bibitem{gottlieb1997gibbs}
D.~Gottlieb and C.-W. Shu.
\newblock On the {Gibbs} phenomenon and its resolution.
\newblock {\em SIAM Review}, 39(4):644--668, 1997.

\bibitem{harten1997high}
A.~Harten.
\newblock High resolution schemes for hyperbolic conservation laws.
\newblock {\em J. Comput. Phys.}, 135(2):260--278, 1997.

\bibitem{harten1987uniformly}
A.~Harten, B.~Engquist, S.~Osher, and S.~R. Chakravarthy.
\newblock Uniformly high order accurate essentially non-oscillatory schemes,
  iii.
\newblock In {\em Upwind and high-resolution schemes}, pages 218--290.
  Springer, 1987.

\bibitem{hesthaven2017numerical}
J.~S. Hesthaven.
\newblock {\em Numerical Methods for Conservation Laws: From Analysis to
  Algorithms}.
\newblock SIAM, 2017.

\bibitem{hesthaven2007nodal}
J.~S. Hesthaven and T.~Warburton.
\newblock {\em Nodal Discontinuous {Galerkin} Methods: Algorithms, Analysis,
  and Applications}.
\newblock Springer Science \& Business Media, 2007.

\bibitem{houston1999posteriori}
P.~Houston, J.~A. Mackenzie, E.~S{\"u}li, and G.~Warnecke.
\newblock A posteriori error analysis for numerical approximations of
  friedrichs systems.
\newblock {\em Numer. Math.}, 82(3):433--470, 1999.

\bibitem{houston2000posteriori}
P.~Houston, R.~Rannacher, and E.~S{\"u}li.
\newblock A posteriori error analysis for stabilised finite element
  approximations of transport problems.
\newblock {\em Comput. Methods Appl. Mech. Eng.}, 190(11-12):1483--1508, 2000.

\bibitem{manteuffel2020least}
D.~Z. Kalchev and T.~A. Manteuffel.
\newblock A least-squares finite element method based on the helmholtz
  decomposition for hyperbolic balance laws.
\newblock {\em arXiv preprint arXiv:1911.05831v2}, 2020.

\bibitem{kingma2015}
D.~P. Kingma and J.~Ba.
\newblock {ADAM}: A method for stochastic optimization.
\newblock In {\em International Conference on Representation Learning, San
  Diego}, 2015; arXiv preprint arXiv:1412.6980.

\bibitem{leveque1992numerical}
R.~J. LeVeque.
\newblock {\em Numerical Methods for Conservation Laws}.
\newblock Birkh\"{a}user, Boston, 1992.

\bibitem{LiuCai1}
M.~Liu, Z.~Cai, and J.~Chen.
\newblock Adaptive two-layer {ReLU} neural network: I. best least-squares
  approximation.
\newblock {\em Computer and Mathematics with Applications}, to appear;
  arXiv:2107.08935v1 [math.NA].

\bibitem{raissi2017physics}
M.~Raissi, P.~Perdikaris, and G.~E. Karniadakis.
\newblock Physics informed deep learning (part {I}): Data-driven solutions of
  nonlinear partial differential equations.
\newblock {\em arXiv preprint arXiv:1711.10561}, 2017.

\bibitem{raissi2019physics}
M.~Raissi, P.~Perdikaris, and G.~E. Karniadakis.
\newblock Physics-informed neural networks: A deep learning framework for
  solving forward and inverse problems involving nonlinear partial differential
  equations.
\newblock {\em J. Comput. Phys.}, 378:686--707, 2019.

\bibitem{roe1981approximate}
P.~L. Roe.
\newblock Approximate riemann solvers, parameter vectors, and difference
  schemes.
\newblock {\em J. Comput. Phys.}, 43(2):357--372, 1981.

\bibitem{shu1998essentially}
C.-W. Shu.
\newblock Essentially non-oscillatory and weighted essentially non-oscillatory
  schemes for hyperbolic conservation laws.
\newblock In {\em Advanced Numerical Approximation of Nonlinear Hyperbolic
  Equations}, pages 325--432. Springer, 1998.

\bibitem{shu1988efficient}
C.-W. Shu and S.~Osher.
\newblock Efficient implementation of essentially non-oscillatory
  shock-capturing schemes.
\newblock {\em J. Comput. Phys.}, 77(2):439--471, 1988.

\bibitem{Sirignano18}
J.~Sirignano and K.~Spiliopoulos.
\newblock {DGM}: A deep learning algorithm for solving partial differential
  equations.
\newblock {\em J. Comput. Phys.}, 375:1139--1364, 2018.

\bibitem{thomas2013numerical}
J.~W. Thomas.
\newblock {\em Numerical Partial Differential Equations: Finite Difference
  Methods}, volume~22.
\newblock Springer Science \& Business Media, 2013.

\bibitem{wang2019learning}
Y.~Wang, Z.~Shen, Z.~Long, and B.~Dong.
\newblock Learning to discretize: solving {1D} scalar conservation laws via
  deep reinforcement learning.
\newblock {\em arXiv preprint arXiv:1905.11079}, 2019.

\end{thebibliography}

\end{document}